\documentclass[12pt]{amsart}

\usepackage{geometry}                
\usepackage{graphicx}
\usepackage{amssymb}

\newtheorem{theorem}{Theorem}[section]

\newtheorem{lemma}[theorem]{Lemma}

\newtheorem*{theoremA}{Theorem A}
\newtheorem*{theoremB}{Theorem B}
\newtheorem*{theoremC}{Theorem C}
\newtheorem*{corollaryD}{Corollary D}
\newtheorem*{corollaryE}{Corollary E}
\newtheorem{corollary}[theorem]{Corollary}
\newtheorem{proposition}[theorem]{Proposition}

\theoremstyle{definition}

\newtheorem{question}[theorem]{Question}

\theoremstyle{remark}
\newtheorem{remark}[theorem]{Remark}

\numberwithin{equation}{section}

\newcommand{\cA}{\mathcal{A}}

\newcommand{\cH}{\mathcal{H}}
\newcommand{\cC}{\mathcal{C}}

\newcommand{\D}{\mathbb{D}}
\newcommand{\B}{\mathbb{B}}
\newcommand{\C}{\mathbb{C}}

\newcommand{\St}{\mathbb{S}}
\newcommand{\R}{\mathbb{R}}

\newcommand{\ga}{\gamma }

\newcommand{\capa}{\operatorname{Cap}}

\newcommand{\bd}[1]{\partial #1}

\newcommand{\Mod}{\operatorname{Mod}}

\newcommand{\card}{\operatorname{card}}

\newcommand{\defeq}{\mathrel{\mathop:}=}

\title[Egg-yolk principle and exponential integrability]{An egg-yolk principle and exponential integrability for quasiregular mappings}

\author{Pietro Poggi-Corradini}
\address{Department of Mathematics, Cardwell Hall, Kansas State University,
Manhattan, KS 66506, USA}
\email{pietro@math.ksu.edu}

\author{Kai Rajala}
\address{Department of Mathematics and Statistics, P.~O.~Box 35 (MaD), FIN-40014, Univ. of Jyv\"askyl\"a, Finland.}
\email{kirajala@maths.jyu.fi}

\thanks{
The second author was supported by the Academy of Finland. Part of this research was done when the second author 
was visiting the University of Cincinnati and the University of
Michigan and the first author was visiting the University of
Michigan. We wish to thank both departments for their hospitality.}
\date{}     
\begin{document}
\begin{abstract}
Quasiregular mappings $f:\Omega\subset\R^{n}\rightarrow \R^{n}$ are a
natural generalization of analytic functions from complex analysis and
provide a theory which is rich with new phenomena. In
this paper we extend a well-known result of A.~Chang and D.~Marshall on
exponential integrability of analytic functions in the disk, to the
case of quasiregular mappings defined in the unit ball of $\R^n$. To
this end, we first establish an ``egg-yolk'' principle for such maps,
which extends a recent result of the first author. Our work leaves
open an interesting problem regarding $n$-harmonic functions.
\end{abstract}

\maketitle
\baselineskip=18pt

\section{Introduction}
We will denote an $n$-dimensional ball with center $a$ and radius $r$ by $\B^n(a,r)$. The unit ball is $\B^n$. Sometimes 
the notation $r\B^n$ for $\B^n(0,r)$ is used. Similarly, 
the notations $\St^{n-1}(a,r)$ and $\St^{n-1}$ for the corresponding $(n-1)$-spheres will be used, respectively. The $s$-dimensional 
Hausdorff measure will be denoted by $\cH_{s}$. The volume of $\B^n$ is denoted by $\alpha_n$, and the $(n-1)$-measure of $\St^{n-1}$ 
by $\omega_{n-1}$.

A mapping $f:\Omega\subset\R^n\rightarrow \R^n$ is called {\sf quasiregular (qr)} if it belongs to
the Sobolev class $W_{loc}^{1,n} (\Omega,\R^n)$, and, for some $K\geq 1$, it satisfies the distortion inequality 
\[
\|Df(x)\|^{n}\leq KJ (x,f)
\]
for almost every $x\in \Omega$, where $\|Df(x)\|$ is the operator norm of the matrix derivative 
$Df(x)=\left(\frac{\partial f_i}{\partial x_j}\right)_{i,j=1}^n$, which is well-defined for almost every $x\in\R^n$, and $J
(x,f)$ is the Jacobian determinant of $f$ at $x$, i.e., $J(x,f)=\det Df(x)$. It is well-known that
quasiregular mappings are continuous and almost everywhere differentiable, and, when
non-constant, they are open and discrete. Also when $n=2$ and $K=1$ they are analytic functions. They provide a
fruitful generalization of classical function theory to higher (real)
dimensional spaces. We refer to \cite{reshetnyak1989} and \cite{rickman1993} for the basic theory of quasiregular mappings. 
The theory of these mappings 
is often referred to, in colorful language, as the {\sf quasiworld}.

The purpose of this paper is twofold. We extend the exponential
integrability result of \cite{chang-marshall:1985ajm} to the
quasiworld. But, to do this, we also need to extend an
``egg-yolk principle for the inverse map'' conjectured by D.~Marshall
in \cite{marshall:1989arkiv}, 
which has been shown to hold in the classical case in \cite{pc:modarg}.

\subsection{Exponential integrability}

The following result is proved in \cite{chang-marshall:1985ajm}.

\begin{theoremA}[Chang-Marshall \cite{chang-marshall:1985ajm}]\label{thm:chang-m}
There is a universal constant $C<\infty$ so that if $f$ is analytic in $\D$, $f(0)=0$, and 
\begin{equation}\label{eq:dirichlet}
\int_\D |f^\prime(z)|^2 \,dA(z)/\pi \leq 1, 
\end{equation}
then 
\[
\int_0^{2\pi}\exp\left(|f^\star(e^{i\theta})|^2\right)\,d\theta \leq C.
\]
where $f^\star$ is the trace of $f$ on $\bd\D$, i.e., $f^\star(\zeta)=\lim_{t\uparrow 1}f(t\zeta)$ for $\cH_1$-a.e.~$\zeta\in \bd\D$.
\end{theoremA}

This result is moreover ``sharp''. Indeed, even though for any given $\beta>0$ and any analytic function $f$ on $\D$, satisfying $f(0)=0$ 
and \eqref{eq:dirichlet}, the integral
\[
\int_0^{2\pi}\exp\left(\beta |f^\star(e^{i\theta})|^2\right)\,d\theta 
\]
is finite, there is a family of functions, the Beurling functions 
\[
B_a(z)=\left(\log\frac{1}{1-az}\right)\left(\log\frac{1}{1-a^2}\right)^{\frac{-1}{2}} \qquad 0<a<1
\]
that are analytic in $\D$, satisfy $B_a(0)=0$ and \eqref{eq:dirichlet}, with the property that for any given $\alpha>1$, 
one can choose $a$ so that the integral
\[
\int_0^{2\pi}\exp\left(\alpha |B_a(e^{i\theta})|^2\right)\,d\theta 
\]
is as large as desired.

In this paper we extend the Chang-Marshall result to quasiregular mappings. 

\begin{theorem}\label{thm:expint}
There exists a constant $C=C(n,K)<\infty$ so that if $f:\B^n\rightarrow\R^n$, $n\geq 2$, is a $K$-quasiregular mapping with $f(0)=0$ and 
\begin{equation}\label{eq:jacob}
\int_{\B^n}J(x,f)\, dx \leq \alpha_n,
\end{equation}
then
\[
\int_{\St^{n-1}}\exp\left((n-1)\Big(\frac{n}{2K}\Big)^{\frac{1}{n-1}}|f^\star(\zeta)|^{\frac{n}{n-1}} \right)\,d\cH_{n-1}(\zeta)\leq C,
\]
where $f^\star$ is the trace of $f$ on $\St^{n-1}$, i.e., $f^\star(\zeta)=\lim_{t\uparrow 1}f(t\zeta)$ for $\cH_{n-1}$-a.e.~$\zeta\in \St^{n-1}$.
\end{theorem}

The trace $f^\star$ in Theorem \ref{thm:expint} is well-defined, since a quasiregular mapping $f:\B^n \rightarrow \R^n$ satisfying \eqref{eq:jacob} 
has radial limits at almost every $\theta \in \St^{n-1}$, see \cite{rickman1993}, VII Theorem 2.7. 

For a mapping satisfying the assumptions of Theorem \ref{thm:expint}, 
\[
\int_{\St^{n-1}}\exp\left(\beta|f^\star(\zeta)|^{\frac{n}{n-1}} \right)\,d\cH_{n-1}(\zeta) < \infty
\]
for every $\beta >0$. This is a consequence of Theorem \ref{thm:bemodest} as will be shown at the end of Section \ref{sec:expint}. 

Theorem \ref{thm:expint} is sharp for $n=2$, in the sense that for any $K \geq 1$ the constant $K^{-1}$ 
cannot be improved on. To see this, first map the unit disk onto the upper half plane by a M\"obius transformation, 
so that $(1,0)$ is mapped to the origin. Then apply the radial stretching $z \mapsto z |z|^{K-1}$, which is a $K$-quasiconformal map, 
and map back to the disk. Finally, apply the Beurling functions $B_a$. The compositions of these maps, $B_{K,a}$, are $K$-quasiregular 
maps satisfying the assumptions of Theorem \ref{thm:expint}, and for each $\beta > K^{-1}$, 
\[
\sup_{0<a<1} \int_0^{2 \pi}\exp\left(\beta|B_{K,a}^\star(e^{i \theta})|^{2} \right)\,d \theta = \infty. 
\]

In dimensions higher than two the situation is different. Indeed, by the Liouville theorem of Gehring and Reshetnyak, see \cite{reshetnyak1989}, Theorem 5.10, 
$1$-quasiregular mappings in dimensions three or higher are M\"obius transformations. Moreover, the $L^{\infty}$-norm of a 
M\"obius transformation satisfying the assumptions of Theorem \ref{thm:expint} is bounded by two. We expect that the constant 
$(n-1)\left(\frac{n}{2K}\right)^{\frac{1}{n-1}}$ is not sharp for any $n \geq 3$ and any $K\geq 1$. In particular, it would be interesting to determine whether the sharp constant stays bounded as $n$ tends to infinity. 
Spatial maps that are similar to the Beurling functions can be constructed by using cylinder maps ($K$-quasiconformal maps mapping $\B^n$ 
onto an infinite cylinder). The best dilatation constant $K$ for cylinder maps is not known, see \cite{gehring-vaisala:acta1965}, Section 8.

\subsection{Further remarks}
The Chang-Marshall theorem has the following two corollaries for harmonic and Sobolev functions.

\begin{corollaryD}\label{cor:harmonic}
There is a universal constant $C<\infty$ so that if $u:\D\rightarrow \R$ is harmonic with $u(0)=0$ and 
\[
\int_\D |\nabla u(z)|^2 \,dA(z)/\pi\leq 1, 
\]
then 
\[
\int_{0}^{2 \pi} \exp\left(u^\star(e^{i \theta})^2\right) \, d\theta \leq C
\]
where $u^\star$ is the trace of $u$ on $\bd\D$, i.e., 
$u^\star(\zeta)=\lim_{t\uparrow 1}u(t\zeta)=u^\star(\zeta)$ for $\cH_1$-a.e. $\zeta\in\bd\D$. 
\end{corollaryD}

\begin{proof}
Let $\tilde{u}$ be the harmonic conjugate of $u$ such that
$\tilde{u}(0)=0$. Then $f=u+i\tilde{u}$ satisfies the hypothesis of
Theorem C, since $|f^\prime|=|\nabla u|$. So
\[
\int_{0}^{2 \pi}\exp\left(u^\star(e^{i \theta})^2\right)\,d\theta \leq 
\int_{0}^{2 \pi}\exp(u^\star(e^{i \theta})^2+\tilde{u}^\star(e^{i \theta})^2)\,d \theta \leq C.
\]
\end{proof}

\begin{corollaryE}\label{cor:sobolev}
There is a universal constant $C<\infty$ so that if $v\in W^{1,2}(\D)$  with $\int_{\bd\D} v^\star(e^{i \theta})\,d\theta=0$ and 
\[
\int_\D |\nabla v(z)|^2 \,dA(z)/\pi\leq 1, 
\]
then 
\[
\int_{0}^{2 \pi}\exp \left(v^\star(e^{i \theta})^2\right) \, d\theta \leq C
\]
where $v^\star$ is the Sobolev trace of $v$ on $\bd\D$.
\end{corollaryE}

For the concept of Sobolev trace see \cite{ziemer1989}, pages 189--191.
\begin{proof}
Let $v^\star$ be the trace of $v$ on the circle $\bd\D$. Solve the
Dirichlet problem with these boundary values, to get $u$ harmonic in
$\D$ with
\[
\int_{\D} |\nabla u|^2 dA/\pi \leq \int_{\D} |\nabla v|^2 \,dA/\pi \le 1.
\]
Then Corollary~D  implies $\int_{0}^{2 \pi} \exp \left(u^\star(e^{i \theta})^2\right)\, d\theta \leq C$, 
but $u^\star=v^\star$. So the same is true for $v^\star$.
\end{proof}

\begin{remark}\label{rem:sobolev}
In terms of statements we have:
\[
\mbox{Theorem~A }\Longrightarrow\mbox{ Corollary~D }
\Longleftrightarrow\mbox{ Corollary~E}
\]
Corollary~E could possibly be proved by ``Sobolev'' methods, see for instance the similar Theorem 3.2.1 of 
\cite{adams-hedberg1996}. When a seemingly stronger normalization 
\[
\int_{\frac{1}{2}\B^n}u(x)\,dx=0 
\]
is assumed, the techniques below can be used to prove results like Corollary~E in all dimensions, see comments 
at the end of Section \ref{sec:capsym}. 
\end{remark}

\begin{remark}\label{rem:essen}
Condition \eqref{eq:dirichlet} says that the Euclidean area of $f(\D)$ counting multiplicity is less or equal to $\pi$. 
In \cite{essen:1987arkiv} it is shown that \eqref{eq:dirichlet} can be replaced by the condition that the area of the set 
$f(\D)$ is less or equal to $\pi$, without counting multiplicity.
\end{remark}

\subsection{Open Questions}
In view of Corollary~D we ask:
\begin{question}
What is the best constant $\beta$ for which there exists $C>0$ so that if $u\in W^{1,n}(\B^n)$, $n\geq 2$, is $n$-harmonic on $\B^n$, $u(0)=0$, and 
\[
\int_{\B^n}|\nabla u(x)|^n\, dx\leq \alpha_n, 
\]
then
\[
\int_{\St^{n-1}}\exp\left(\beta |u^\star(\zeta)|^{\frac{n}{n-1}} \right)\,d\cH_{n-1}(\zeta)\leq C ?
\]
\end{question}

\subsection{Beurling's estimate}
In \cite{marshall:1989arkiv}, Don Marshall deduces Theorem~A from an estimate of Beurling, Theorem~B below.
We denote $E_t=\{x \in \B^n: |f(x)|=t\}$, and $F^\star_s=\{\theta \in \St^{n-1}: |f(\theta)|>s\}$. 
The following is an unpublished estimate of A.~Beurling which  is stated and proved in 
\cite{marshall:1989arkiv}. Here ``$\capa$'' denotes logarithmic capacity.

\begin{theoremB}[Beurling]\label{thm:beurling}
Suppose $f$ is analytic in a neighborhood of $\overline{\D}$ and suppose that $|f(z)|\leq M$ for $|z|\leq r<1$, for some $0<r<1$. Then, for every $s>M$,
\[\capa F^\star_s\leq r^{\frac{-1}{2}}\exp\left(-\pi\int_M^s\frac{dt}{|f(E_t)|}\right)\]
where $|f(E_t)|$ denotes the length of $f(E_t)$ counting multiplicity.
\end{theoremB}

We establish a similar estimate in space. For a quasiregular map $f:\B^n \rightarrow \R^n$, $n \geq 2$, we denote the 
$(n-1)$-measure of $f(E_t)$ counting multiplicity by $\cA_{n-1}f(E_t)$; 
$$
\cA_{n-1}f(E_t)= \int_{\St^{n-1}(0,t)}\card f^{-1}(y)\, d\cH_{n-1}(y). 
$$

\begin{theorem}\label{thm:bemodest}
Let $f$ be a $K$-quasiregular mapping defined in a neighborhood of $\overline{\B}^n$, $n \geq 2$, and suppose that $|f(x)|\leq M$ for $|x|\leq r <1$. 
Then, for every $s>M$,
\begin{equation}\label{eq:expdec}
\cH_{n-1}(F^\star_s)\leq C_1\exp\left( (1-n)\Big(\frac{\omega_{n-1}}{2K} \Big)^{\frac{1}{n-1}}\int_M^s\frac{dt}{(\cA_{n-1}f(E_t))^{\frac{1}{n-1}}}\right), 
\end{equation}
where $C_1$ depends only on $n$, $K$ and $r$.
\end{theorem}

\subsection{An egg-yolk principle for the inverse}
In \cite{marshall:1989arkiv}, Don Marshall conjectures an egg-yolk
principle that would have simplified his argument for passing from Theorem~B to Theorem~A.  
This was proved in \cite{pc:modarg} by the first author.
\begin{theoremC}[\cite{pc:modarg}]
There is a universal constant $0<r_0<1$ such that whenever $f$ is analytic on
$\D\defeq\{z\in \C: |z|<1\}$ with $f(0)=0$, and whenever $M>0$ is such that 
\[
\int_{\{z\in \D:|f(z)|<M\}}|f^\prime(z)|^2\,dA(z) <\pi M^2, 
\]
then we have that $|z|<r_0$ implies $|f(z)|<M$.
\end{theoremC}

Here we prove that Theorem~C extends to quasiregular maps, and this
will allow us to deduce Theorem \ref{thm:expint} from Theorem \ref{thm:bemodest}.

\begin{theorem}\label{thm:kai}
Given $n\geq 2$ and $K\geq 1$,
there exists a constant $0<r_0(n,K)<1$, so that whenever $f:\B^n\rightarrow\R^n$ is a $K$-quasiregular 
mapping with $f(0)=0$ and whenever $M>0$ is such that
\begin{equation}
\label{eq:mult}
\int_{\{x\in\B^n: |f(x)|<M \}}J(x,f)\, dx < \alpha_n M^n, 
\end{equation}
then we have that $|x|<r_0$ implies $|f(x)|<M$.
\end{theorem}

Theorem \ref{thm:kai} is equivalent to the following.

\begin{corollary}\label{cor:equiv}
For $n\geq 2$ and $K\geq 1$,
there exists a constant $0<r_0(n,K)<1$ so that if $f:\B^n\rightarrow\R^n$ is a $K$-quasiregular 
mapping with $f(0)=0$, then 
$0\leq M<\max_{|x|\leq r_0}|f(x)|$ implies 
\[
\int_{\{x\in \B^n:|f(x)|<M\}}J(x,f)\,dx \geq  \alpha_n M^n. 
\]
\end{corollary}

Theorem \ref{thm:kai} no longer holds true if instead of \eqref{eq:mult} it is assumed that $\B^n\setminus f(\B^n) \neq \emptyset$, see 
\cite{pc:modarg}, Remark 1.5. 

\vskip 15pt

\noindent{\bf{Acknowledgements}}

\noindent We thank Pekka Koskela for useful discussions.

\section{Proof of Theorem \ref{thm:kai}}

We first recall the classical (conformal) modulus for path families in $\R^n$.
Let $\Gamma$ be a family of paths $\ga$, i.e., continuous functions $\ga:I \rightarrow \R^n$, where $I=[a,b]$ or $[a,b)$. We say that a Borel measurable 
function $\rho:\R^n\rightarrow [0,+\infty]$ is {\sf admissible for $\Gamma$} if
\[
\int_\ga \rho \, ds \geq 1\qquad \forall \ga\in\Gamma.
\]
Then the {\sf modulus} of $\Gamma$ is
\[
\Mod \Gamma\defeq\inf\left\{\int_{\R^n}\rho(x)^n\,dx : \rho\mbox{ admissible} \right\}.
\]
%

We recall two classical results concerning conformal modulus.

\begin{lemma}[Poletsky's inequality, \cite{rickman1993}, II Theorem 8.1]
\label{pol}
Let $f: \Omega \to \R^n$ be a non-constant $K$-quasiregular mapping, and $\Gamma$ a family of paths 
in $\Omega$. Then 
$$
\Mod f\Gamma \leq K^{n-1} \Mod \Gamma. 
$$
\end{lemma}

\begin{lemma}[\cite{vaisala1971}, Theorem 10.12]
\label{geh}
Suppose that $J$ is a measurable set of radii, and $p \in \R^n$. For each $r \in J$, consider distinct points 
$a_r, b_r$ in $\St^{n-1}(p,r)$. Set 
$$
\Gamma = \{ \gamma : [a,b) \to \St^{n-1}(p,r)|\, r \in J, \, \gamma \text{ connects } a_r \text{ and } b_r \}. 
$$
Then 
$$
\Mod \Gamma \geq c_n \int_J \frac{dr}{r},  
$$
where $c_n >0$ only depends on $n$. 
\end{lemma}

Let $f$ satisfy the assumptions of Theorem \ref{thm:kai}. We lose no generality by assuming $M=1$. 
Let $\delta$ denote the largest radius so that 
$$
f(\B^n(0,\delta)) \subset \B^n. 
$$
In order to prove Theorem \ref{thm:kai} we need to show that $\delta \geq C(n,K)$. Also, we let
\begin{eqnarray*}
F_0 &=& \B^n \setminus f(\B^n), \\
F_1 &=& \{ y \in \B^n: \card f^{-1}(y)=1 \}, \\ 
F_m &=& \{ y \in \B^n: \card f^{-1}(y) \geq 2 \}= \B^n \setminus (F_0 \cup F_1).  
\end{eqnarray*}
By (\ref{eq:mult}) and a change of variables, we have
$$
\alpha_n > \int_{\{x \in \B^n: f(x) \in \B^n\}} J(x,f)\, dx = \int_{\B^n} \card f^{-1}(y)\,dy. 
$$
Therefore $F_0 \neq \emptyset$. 

We first prove Theorem \ref{thm:kai} under the assumption 
\begin{equation}
\label{case1}
|F_0| \geq \alpha_n 100^{-n}.
\end{equation}
We denote by $T$ the set of those radii $0<r<1$ for which 
$$
\St^{n-1}(0,r) \cap F_0 \neq \emptyset.
$$
 
\begin{lemma}
\label{polar}
Assume that \eqref{case1} holds true. Then 
$$
\int_T \frac{dr}{r} \geq n^{-1} 100^{-n}. 
$$
\end{lemma}
\begin{proof}
Since $r<1$, we have 
\begin{eqnarray*}
\int_T \frac{dr}{r} &=& \omega_{n-1}^{-1} \int_T \int_{\St^{n-1}(0,r)} r^{-n} \,d\cH_{n-1} \, dr 
\geq \omega_{n-1}^{-1} \int_{\R^n} \chi_{\{y: |y| \in T\}}(x)\,dx \\
&=& \omega_{n-1}^{-1} |\{y: |y| \in T\}| \geq \omega_{n-1}^{-1} |F_0| 
\geq \alpha_n \omega_{n-1}^{-1} 100^{-n} = n^{-1}100^{-n}.
\end{eqnarray*}

\end{proof}

\begin{proposition}
\label{pcase1} 
Theorem \ref{thm:kai} holds true under assumption \eqref{case1}. 
\end{proposition}

\begin{proof}
By definition of $T$, for each $r \in T$, we can choose points $q_r \in F_0\cap \St^{n-1}(0,r)$. Also, since $\overline{f(\B^n(0,\delta))}$ is a connected set containing 
$0$ and a point in $\St^{n-1}$, for each $r\in T$, we can choose points $a_r\in \B^n(0,\delta)$ such that $f(a_r)\in \St^{n-1}(0,r)$.
Then, for every path $\gamma$ starting at $f(a_r)$ and joining 
$f(a_r)$ to $q_r$ in $\St^{n-1}(0,r)$, every maximal lift $\gamma'$ of $\gamma$ starting at $a_r$ 
accumulates on $\St^{n-1}$ (see \cite{rickman1993}, II.3 for the definition of a maximal lift). Hence, if we denote the family of all 
such lifts, for any $r \in T$, by $\Gamma$, we have 
\begin{equation}
\label{u1}
\Mod \Gamma \leq \omega_{n-1} \left(\log \delta^{-1}\right)^{1-n}. 
\end{equation} 
On the other hand, by Lemmas \ref{geh} and \ref{polar}, 
\begin{equation}
\label{l1}
\Mod f \Gamma \geq c_n \int_T \frac{dr}{r} \geq c_n n^{-1} 100^{-n}. 
\end{equation}
By combining \eqref{u1}, \eqref{l1} and Lemma \ref{pol}, we have 
$$
c_n n^{-1} 100^{-n} \leq K^{n-1} \omega_{n-1} \left(\log \delta^{-1}\right)^{1-n}, 
$$
Thus Theorem \ref{thm:kai} holds in this case with
$$
r_0(n,K)= \exp \Big( -\Big(100^n c_n^{-1}nK^{n-1}\omega_{n-1} \Big)^{\frac{1}{n-1}} \Big). 
$$
\end{proof}

We now treat the case when (\ref{case1}) fails. First we establish a geometric lemma. 

\begin{lemma}
\label{geom}
Fix $q \in F_0$. Then there exists a point $w \in \B^n$, and $1/4  \leq s < 1$, such that for all $r \in (s, \sqrt{3}s)$, 
we have $q \in \B^n(w,r)$ and $\St^{n-1}(w,r)\cap f(\B^n(0,\delta))\neq \emptyset$.
\end{lemma}

\begin{proof}
First assume $|q| \leq 1/2$. Then, since $\overline{f(\B^n(0,\delta))}$ is a connected set containing 
$0$ and a point in $\St^{n-1}$, 
$$
\St^{n-1}(0,r) \cap f(\B^n(0,\delta)) \neq \emptyset \quad \forall \, r \in \left(\frac{1}{2}, \frac{\sqrt{3}}{2}\right). 
$$
Hence we may choose $w=0$, $s=1/2$. 

Thus assume $|q| > 1/2$. Choose $p \in \B^n(0,\delta)$ such that $|f(p)|=|q|$. Consider the triangle 
with vertices $0$, $f(p)$ and $q/2$. Then, if the angle of the triangle at $q/2$ is less than $\pi /2$, 
we have, for each $r \in (|q|/2, \sqrt{3}|q|/2)$, 
$$
0,q \in \B^n(q/2,r), \quad f(p) \notin \B^n(q/2,r). 
$$
Since $0,f(p) \in f(\B^n(0,\delta))$, there exists, for each such $r$, a point $\eta_r \in \B^n(0,\delta)$ 
such that $f(\eta_r) \in \St^{n-1}(q/2,r)$. Hence we may choose $w=q/2$ and $s= |q|/2 > 1/4$ in this case. 

If the angle is greater than or equal to $\pi/2$, we have, for each $r \in (|q|/2, \sqrt{3}|q|/2)$,  
$$
f(p), q \in \B^n\left(\frac{f(p)+q}{2},r\right), \quad 0 \notin \B^n\left(\frac{f(p)+q}{2},r\right). 
$$ 
Hence we may in this case choose $w= (f(p)+q)/2$ and $s=|q|/2$. 
\end{proof}

Let $q$, $w$, and $s$ be as in Lemma \ref{geom}. We denote by $G$ the set of all radii $r\in (s, \sqrt{3}s)$ 
for which 
$$
F_1 \cap \St^{n-1}(w,r) \neq \emptyset. 
$$ 

\begin{lemma}
\label{polar2}
If \eqref{case1} fails, then 
$$
\int_G \frac{dr}{r} \geq n^{-1}100^{-n}. 
$$
\end{lemma}

\begin{proof}
As in Lemma \ref{polar}, we have 
\begin{equation}
\label{tsup}
\int_G \frac{dr}{r} \geq \omega_{n-1}^{-1} \int_G \int_{\St^{n-1}(w,r)}\, d\cH_{n-1} \, r^{-n}\,dr 
\geq \omega_{n-1}^{-1} |F_1 \cap (\B^n(w,\sqrt{3}s)\setminus \overline{\B^n}(w,s))|. 
\end{equation}
By our assumption (\ref{eq:mult}) and a change of variables, 
$$
|F_1|+ 2|F_m| \leq \int_{\B^n} \card f^{-1}(y)\,dy = \int_{\{x \in \B^n: f(x) \in \B^n\}} J(x,f)\,dx < \alpha_n 
= |F_0|+|F_1|+|F_m|.  
$$
So
\begin{equation}
\label{so}
|F_m| \leq |F_0| \leq \alpha_n 100^{-n}, 
\end{equation}
where the last inequality holds true since we assume the converse of \eqref{case1}. 

On the other hand, since $w \in \B^n$ and $s \geq 1/4$, we have 
\begin{equation}
\label{su}
|(\B^n(w,\sqrt{3}s)\setminus \overline{\B^n}(w,s)) \cap \B^n| \geq \alpha_n10^{-n}, 
\end{equation}
and combining \eqref{so} and \eqref{su} yields 
\begin{eqnarray} 
\label{si}
& &|F_1 \cap (\B^n(w,\sqrt{3}s)\setminus \overline{\B^n}(w,s))|=
|(\B^n(w,\sqrt{3}s)\setminus \overline{\B^n}(w,s))\cap \B^n| \\ \nonumber &-& 
|(\B^n(w,\sqrt{3}s)\setminus \overline{\B^n}(w,s))\cap (F_0 \cup F_m)| \geq \alpha_n 100^{-n}. 
\end{eqnarray}
The Lemma follows by combining \eqref{tsup} and \eqref{si}. 
\end{proof}

For each $r \in G$, choose points $p_r \in f^{-1}(F_1)$, $a_r \in \B^n(0,\delta)$ such that 
$$
f(p_r), f(a_r) \in \St^{n-1}(w,r).
$$ 
Denote 
\begin{eqnarray*}
G_1 &=& \{r \in G: |p_r| \geq \delta^{\frac{1}{2}} \}, \\
G_2 &=& \{r \in G: |p_r| < \delta^{\frac{1}{2}}\} = G \setminus G_2.  
\end{eqnarray*} 
Then, by Lemma \ref{polar2}, either \eqref{case1} holds, or else we have one of 
\begin{eqnarray}
\label{case2}
\int_{G_1} \frac{dr}{r} &\geq& 2^{-1}n^{-1}100^{-n}, 
\end{eqnarray}
or 
\begin{eqnarray}
\label{case3}
\int_{G_2} \frac{dr}{r} &\geq& 2^{-1}n^{-1}100^{-n}.  
\end{eqnarray}

\begin{proposition}
\label{pcase2}
Theorem \ref{thm:kai} holds true under assumption \eqref{case2}. 
\end{proposition}

\begin{proof}
For each $r \in G_1$ and each $\gamma$ starting at $f(a_r)$ and joining $f(a_r)$ to $f(p_r)$ in 
$\St^{n-1}(w,r)$, consider a maximal lift $\gamma'$ of $\gamma$ starting at $a_r$. Then, since 
$\card f^{-1}(f(p_r))=1$, either $\gamma'$ accumulates to $\St^{n-1}$, or $\gamma'$ ends at $p_r$; in any case, 
$\gamma'$ starts at $\B^n(0,\delta)$ and leaves $\B^n(0,\delta^{\frac{1}{2}})$. Denote the family of all such 
$\gamma'$ by $\Gamma$. Then we have 
\begin{equation}
\label{u2}
\Mod \Gamma \leq \omega_{n-1} \left(\log \frac{\delta^{\frac{1}{2}}}{\delta}\right)^{1-n} = \omega_{n-1} \left(\log \delta^{\frac{-1}{2}}\right)^{1-n}. 
\end{equation}

On the other hand, combining Lemma \ref{geh} and \eqref{case2} yields 
\begin{equation}
\label{l2}
\Mod f \Gamma \geq c_n 2^{-1}n^{-1}100^{-n}. 
\end{equation}
Furthermore, combining \eqref{u2}, \eqref{l2} and Lemma \ref{pol} gives 
$$
c_n 2^{-1}n^{-1}100^{-n} \leq K^{n-1} \omega_{n-1} \left( \log \delta^{\frac{-1}{2}}\right)^{1-n}, 
$$
Thus Theorem \ref{thm:kai} holds in this case with
$$
r_0(n,K)= \exp \Big( -2\Big(100^n 2c_n^{-1}nK^{n-1}\omega_{n-1} \Big)^{\frac{1}{n-1}} \Big). 
$$
\end{proof}

In order to finish the proof of Theorem \ref{thm:kai}, we need the following auxiliary result. 
\begin{lemma}
\label{to}
For each $r \in G_2$ there exists $\tau_r \in \St^{n-1}(0,\delta^{\frac{1}{4}})$ such that $f(\tau_r)\in \St^{n-1}(w,r)$. 
\end{lemma}

\begin{proof}
Let $U_r$ be any component of $f^{-1}(\B^n(w,r))$ intersecting $\B^n(0,\delta)$. Such a component exists by 
Lemma \ref{geom}. Also, by Lemma \ref{geom}, $\B^n(w,r) \setminus f(\B^n)\neq \emptyset$, and hence 
$f_{|U_r}:U_r \to \B^n(w,r)$ is not onto. Thus, by \cite{rickman1993}, I Lemma 4.7, 
$$
\St^{n-1}(0,t) \cap U_r \neq \emptyset \quad \forall t \in (\delta,1). 
$$

Choose $k_r \in U_r \cap \St^{n-1}(0,\delta^{\frac{1}{4}})$, and consider all paths joining $k_r$ to $-k_r$ in 
$\St^{n-1}(0,\delta^{\frac{1}{4}})$. If none of the images of these paths intersects $\St^{n-1}(w,r)$, we have 
\begin{equation}
\label{op}
f(\St^{n-1}(0,\delta^{\frac{1}{4}})) \subset \B^n(w,r). 
\end{equation}
Since $f$ is open, 
$$
\partial f(\B^n(0,\delta^{\frac{1}{4}})) \subset f(\St^{n-1}(0,\delta^{\frac{1}{4}})), 
$$
and since $f(\B^n(0,\delta^{\frac{1}{4}}))$ is bounded, \eqref{op} further implies 
\begin{equation}
\label{ope}
f(\B^n(0,\delta^{\frac{1}{4}})) \subset \B^n(w,r). 
\end{equation}
By Lemma \ref{geom} there are, however, points $x \in \B^n(0,\delta)$ such that $f(x)\notin \B^n(w,r)$ 
which contradicts \eqref{ope}. The proof is complete.
\end{proof}

\begin{proposition}
\label{pcase3}
Theorem \ref{thm:kai} holds true under assumption \eqref{case3}. 
\end{proposition}

\begin{proof}
For each $r \in G_2$, and each $\gamma$ starting at $f(\tau_r)$ (where $\tau_r$ is as in 
Lemma \ref{to}) and joining $f(\tau_r)$ to $f(p_r)$ in $\St^{n-1}(w,r)$, consider a maximal lift $\gamma'$ 
of $\gamma$ starting at $\tau_r$. Then, since $\card f^{-1}(f(p_r))=1$, either $\gamma'$ accumulates to $\St^{n-1}$, 
or $\gamma'$ ends at $p_r$. We denote the family of all such $\gamma'$ for which the first case occurs 
by $\Gamma_1$, the family of all $\gamma'$ for which the second case occurs by $\Gamma_2$, and 
$\Gamma = \Gamma_1 \cup \Gamma_2$. 

Then, since each $\gamma' \in \Gamma_1$ connects $\St^{n-1}(0,\delta^{\frac{1}{4}})$ to $\St^{n-1}$, 
\begin{equation}
\label{u31}
\Mod \Gamma_1 \leq \omega_{n-1} \left(\log \delta^{\frac{-1}{4}}\right)^{1-n}. 
\end{equation} 
Similarly, since $p_r \in \B^n(0,\delta^{\frac{1}{2}})$ for all $r \in G_2$, 
\begin{equation}
\label{u32}
\Mod \Gamma_2 \leq \omega_{n-1} \left(\log \frac{\delta^{\frac{1}{4}}}{\delta^{\frac{1}{2}}} \right)^{1-n}
= \omega_{n-1} \left(\log \delta^{\frac{-1}{4}}\right)^{1-n}. 
\end{equation}

By Lemma \ref{geh} and \eqref{case3}, 
\begin{equation}
\label{l3}
\Mod f \Gamma \geq c_n 2^{-1} n^{-1}100^{-n}. 
\end{equation}
Hence, combining \eqref{u31}, \eqref{u32}, \eqref{l3} and Lemma \ref{pol} yields 
\begin{eqnarray*}
c_n 2^{-1}n^{-1}100^{-n} &\leq& \Mod f \Gamma \leq K^{n-1} \Mod \Gamma \leq K^{n-1}(\Mod \Gamma_1+ \Mod \Gamma_2)\\
&\leq & 2K^{n-1} \omega_{n-1} \left(\log \delta^{\frac{-1}{4}}\right)^{1-n},  
\end{eqnarray*}
Thus Theorem \ref{thm:kai} holds in this case with
$$
r_0(n,K)= \exp \Big( -4\Big(100^n 4c_n^{-1}nK^{n-1}\omega_{n-1} \Big)^{\frac{1}{n-1}} \Big). 
$$
\end{proof}

The proof of Theorem \ref{thm:kai} follows by combining Propositions \ref{pcase1}, \ref{pcase2} and \ref{pcase3}.

\section{Beurling's modulus estimate}\label{sec:bemodest}

Suppose $f$ is $K$-quasiregular in a neighborhood of the closed unit ball $\overline{\B}^n$, and for some fixed $0<r<1$ let
 $M\defeq\max_{r\overline{\B}^n}|f|$. Recall that 
for $s>M$ we define $F^\star_s=\{\zeta\in \St^{n-1}: |f(\zeta)|\geq s\}$ and for $M<t<s$ we have  $E_t=\{x\in\B^n: |f(x)|=t\}$. 
Consider the  family $\Gamma_s$ consisting of the paths in $\B^n$ starting at $r\overline{\B}^n$ and ending at 
$F^\star_s$. We claim that
\begin{equation}\label{eq:bemodest}
\Mod \Gamma_s\leq  K\left(\int_M^s\frac{dt}{(\cA_{n-1}f(E_t))^{\frac{1}{n-1}}}\right)^{1-n}. 
\end{equation}
Recall that $\cA_{n-1}f(E_t)= \int_{\St^{n-1}(0,t)} \operatorname{card} f^{-1}(y)\,d\cH_{n-1}(y)$.

\begin{proof}
Set $\rho : \R^n \rightarrow [0, \infty)$, 
\[
\rho(x)= \left(\int_M^s \frac{du}{(\cA_{n-1}f(E_u))^{\frac{1}{n-1}}}\right)^{-1} 
\frac{\|Df(x)\|}{(\cA_{n-1}f(E_t))^{\frac{1}{n-1}}} \quad \text{when }|f(x)|=t \in (M,s), 
\]
and $\rho(x)=0$ otherwise. Then, for each $\gamma \in \Gamma_s$, 
\[
\int_{\gamma} \rho \, ds \geq \left(\int_M^s \frac{du}{(\cA_{n-1}f(E_u))^{\frac{1}{n-1}}}\right)^{-1} 
\int_{f(\gamma)} (\cA_{n-1}f(E_{|\cdot|}))^{\frac{-1}{n-1}}\, ds \geq 1.
\]
Moreover, if we denote 
\[
I(M,s)= \int_M^s \frac{du}{(\cA_{n-1}f(E_u))^{\frac{1}{n-1}}} 
\]
and 
\[
A(M,s)= f^{-1}(\B^n(0,s)\setminus \overline{\B^n}(0,M)), 
\]
we have 
\begin{eqnarray*}
\Mod \Gamma_s &\leq& \int_{\R^n} \rho(x)^n\,dx = I(M,s)^{-n} \int_{A(M,s)} \frac{\|Df(x)\|^n}{(\cA_{n-1}f(E_{|f(x)|}))^{\frac{n}{n-1}}}\,dx  \\ 
&\leq& K I(M,s)^{-n} \int_{A(M,s)} \frac{J(x,f)}{(\cA_{n-1}f(E_{|f(x)|}))^{\frac{n}{n-1}}}\,dx  \\ 
&=& K I(M,s)^{-n} \int_{f(A(M,s))} \frac{\card f^{-1}(y)}{(\cA_{n-1}f(E_{|y|}))^{\frac{n}{n-1}}}\,dy \\
&=& K I(M,s)^{-n} \int_M^s (\cA_{n-1}f(E_{t}))^{\frac{-n}{n-1}} \int_{\St^{n-1}(0,t)} \card f^{-1}(\varphi)\, d \cH_{n-1}(\varphi)\,dt = K I(M,S)^{1-n}.
\end{eqnarray*}

\end{proof}

\section{Capacity and Symmetrization}
\label{sec:capsym}
We recall that a {\sf condenser} is a pair $(\Omega,K)$ with $\Omega\subset \R^n$, $\Omega$ open and $K$ compact with 
$\emptyset\neq K\subset \Omega$. Also, the {\sf conformal capacity} of $(\Omega,K)$ is
\[
\capa(\Omega,K)\defeq\inf\{
\|\nabla u\|_{L^n(\Omega)}^n: u\in W^{1,p}_0(\Omega), u_{|V}\geq 1, \text{ for some } V \mbox{ open }, V\supset K\}
\]
where $W^{1,p}_0(\Omega)$ is the closure of $C^\infty_0(\Omega)$ (the smooth functions compactly supported in $\Omega$) in the norm
\[
\|u\|_{W^{1,p}_0(\Omega)}=\left(\int_\Omega |u(x)|^n+|\nabla u(x)|^n \,dx\right)^{\frac{1}{n}}.
\]
By Proposition II.10.2 of \cite{rickman1993}, if $\Gamma(\Omega,K)$ is the family of all paths 
$\gamma: [a,b)\rightarrow\Omega$ such that $\gamma(a)\in K$ and $\lim_{t\rightarrow b}\gamma(t)\in \partial\Omega$, then
\begin{equation}\label{eq:capmod}
\capa(\Omega,K)=\Mod\Gamma(\Omega,K). 
\end{equation}

We are mainly interested in measuring the sets $F_s^\star$ defined in Section \ref{sec:bemodest}, which are compact subsets of 
$\St^{n-1}$. Therefore, we will fix $0<r<1$ to be determined later, consider the spherical ring $A(r)=\{x\in\R^n: r<|x|<1/r\}$, $0<r<1$, 
and compute $\capa(A(r),F^\star_s)$.

By the symmetry rule, cf. \cite{garnett-marshall2005} IV(3.4), if $F\subset\St^{n-1}$, we have:
\begin{equation}\label{eq:symmetry}
\Mod \Gamma_s =\frac{1}{2}\Mod \Gamma(A(r),F)=\frac{1}{2}\capa(A(r),F).
\end{equation}

Also, if $F\subset\St^{n-1}$, let $\cC(F)$ be the spherical cap centered at $e_1=(1,0,\dots,0)$ with 
$\cH_{n-1}(\cC(F))=\cH_{n-1}(F)$. By spherical symmetrization, see \cite{gehring:1961tams}, 
\begin{equation}\label{eq:sphsym}
\capa(A(r),\cC(F))\leq\capa(A(r),F).
\end{equation}
By \cite{gehring:1961tams}, Theorem 4, we see that, when $\cH_{n-1}(F)\leq \epsilon(r,n)$, 
\begin{equation} \label{eq:another}
\capa(A(r),C(F)) \geq \omega_{n-1} \log^{1-n} \frac{C_2}{\cH_{n-1}(F)^{\frac{1}{n-1}}}, 
\end{equation}
where $C_2>0$ depends only on $r$ and $n$ (the results in \cite{gehring:1961tams} are stated for $n=3$ only, but they extend to all dimensions). 

Putting \eqref{eq:bemodest}, \eqref{eq:capmod}, \eqref{eq:symmetry}, \eqref{eq:sphsym}, and \eqref{eq:another} together, we obtain
\eqref{eq:expdec} and thus we have proved Theorem \ref{thm:bemodest} for 
$\cH_{n-1}(F^\star_s)\leq \epsilon(r,n)$. If $\cH_{n-1}(F^\star_s)> \epsilon(r,n)$, then the arguments above show that 
\begin{equation}
\label{eq:lower}
\Mod \Gamma_s \geq C(r,n). 
\end{equation}
Combining \eqref{eq:lower} with \eqref{eq:bemodest} yields 
\[
\int_M^s\frac{dt}{(\cA_{n-1}f(E_t))^{\frac{1}{n-1}}} \leq C(r,n,K). 
\]
Hence increasing $C_1$ if necessary gives Theorem \ref{thm:bemodest} for all $s>M$.

We finish this section by briefly commenting on the real-valued case mentioned in the introduction. Suppose that 
$u:\B^n \rightarrow \overline{\R}$ belongs to $W^{1,n}(\B^n)$ and satisfies 
\begin{equation}
\label{eq:null}
\int_{\frac{1}{2}\B^n} u(x)\,dx = 0. 
\end{equation}
Then, by the Poincar\'e inequality and \eqref{eq:null}, 
\[
|A_T|=|\{x \in \frac{1}{2}\B^n: |u|\leq T\}| \geq C(n) 
\]
for large enough $T$ depending only on $n$ and the Sobolev norm of $u$. Hence, by applying arguments similar to the ones above to the 
$n$-capacity related to the sets $A_T$ and $U^\star_s=\{y \in \St^{n-1}: |u^\star (y)|\geq s\}$, we have an estimate for the $\cH_{n-1}$-measure of 
$U^\star_s$ in terms of $s$, $T$ and $\int_{\{x \in \B^n:|u|\leq s\}} |\nabla u(x)|^n\,dx$.

\section{Exponential integrability}
\label{sec:expint}

In this section we prove Theorem \ref{thm:expint} by using the results established in previous sections and 
arguments similar to those used in \cite{marshall:1989arkiv}. Let $f$ be a $K$-quasiregular mapping defined in a neighborhood 
of $\B^n$ and satisfying $f(0)=0$ and \eqref{eq:jacob}. We denote 
$$
\beta= (n-1)\Big(\frac{n}{2K} \Big)^{\frac{1}{n-1}}. 
$$
Then 
$$
\alpha_n^{\frac{1}{n-1}}\beta = (n-1)\Big(\frac{\omega_{n-1}}{2K}\Big)^{\frac{1}{n-1}}. 
$$
Notice that we lose no generality by assuming that $f$ is defined in a neighborhood of $\B^n$: if we 
consider a sequence $(r_j)$ increasing to one, and functions $f_j$, $f_j(x)=f(r_jx)$, then the existence of 
radial limits at almost every $\varphi \in \St^{n-1}$ and Fatou's lemma yield 
\[
\int_{\St^{n-1}}\exp\left(\beta |f^\star(\zeta)|^{\frac{n}{n-1}} \right)\,d\cH_{n-1}(\zeta) 
\leq \liminf_j \int_{\St^{n-1}}\exp\left(\beta |f_j^\star(\zeta)|^{\frac{n}{n-1}} \right)\,d\cH_{n-1}(\zeta). 
\]

By the Cavalieri principle, 
\begin{equation}
\label{eq:cava}
\int_{\St^{n-1}} \exp \Big(\beta |f(\zeta)|^{\frac{n}{n-1}}\Big)\, d\cH_{n-1} (\zeta) = \omega_{n-1} + 
\frac{\beta n}{n-1} \int_0^{\infty} s^{\frac{1}{n-1}}\cH_{n-1}(F^\star_s)\exp (\beta s^{\frac{n}{n-1}})\,ds. 
\end{equation}
We choose $r_0=r_0(n,K)$ as in Theorem \ref{thm:kai}, and let $M=\max_{|x|\leq r_0}|f(x)|$. Note that by 
Corollary \ref{cor:equiv} and \eqref{eq:jacob}, we have $M<1$ and 
\begin{equation}
\label{eq:MM}
\int_{\{x \in \B^n:f(x) \in \B^n(0,M)\} } J(x,f)\,dx = \int_0^M \cA_{n-1}f(E_t) \,dt \geq \alpha_n M^n. 
\end{equation}

Using \eqref{eq:expdec} and \eqref{eq:cava}, we are reduced to estimate 
\begin{equation}
\label{eq:redu}
\int_0^{\|f\|_{\infty}} s^{\frac{1}{n-1}}\exp (\beta s^{\frac{n}{n-1}}-\psi(s))\,ds, 
\end{equation}
where $\psi(s)=0$ for $0<s \leq M$ and 
\[
\psi(s)= \alpha_n^{\frac{1}{n-1}} \beta \int_M^s \frac{dt}{(\cA_{n-1}f(E_t))^{\frac{1}{n-1}}} 
\]
for $s \geq M$. We modify $\psi$ as follows: for $0<s \leq M$, set 
\[
\tilde{\psi}(s)\defeq \mu s, 
\]
and for $s \geq M$, 
\[
\tilde{\psi}(s)\defeq \psi(s)+ \mu M, 
\]
where 
\[
\mu = \Big(\frac{M \beta^{n-1} \alpha_n}{\int_0^M \cA_{n-1}f(E_t)\,dt} \Big)^{\frac{1}{n-1}}. 
\]
Note that $\tilde{\psi}$ is strictly increasing for $0<s \leq \|f\|_{\infty}$ and constant, equal to $\|\tilde{\psi}\|_{\infty}$, 
for $s > \|f\|_{\infty}$. Also $\tilde{\psi}(0)=0$. Finally, $\tilde{\psi} \leq \psi + \mu M$. So, by \eqref{eq:MM}, and since $M<1$, 
it is enough to estimate \eqref{eq:redu} with $\psi$ replaced by $\tilde{\psi}$. 

Let $\phi(y)\defeq\tilde{\psi}^{-1}(y)$ for $0<y \leq \|\tilde{\psi}\|_{\infty}$ and $\phi(y)\defeq \|f\|_{\infty}$ for 
$y > \|\tilde{\psi}\|_{\infty}$, so that $\phi$ is strictly increasing for $0<y \leq \|\tilde{\psi}\|_{\infty}$ and $\phi(0)=0$. 

Changing variables $y=\tilde{\psi}(s)$ the integral \eqref{eq:redu} becomes 
\[
\int_0^{\|\tilde{\psi}\|_{\infty}} \exp(\beta \phi(y)^{\frac{n}{n-1}}-y)\phi'(y) \phi(y)^{\frac{1}{n-1}}\,dy 
\]
which, since $\phi'\geq 0$, is less than or equal to the same integral but from $0$ to $\infty$. Integrating by 
parts we then need to estimate 
\begin{equation}
\label{eq:inparts}
\int_0^{\infty} \exp(\beta \phi(y)^{\frac{n}{n-1}}-y)\,dy = \int_0^{\infty} \exp( (\beta^{\frac{n-1}{n}} \phi(y))^{\frac{n}{n-1}}-y)\,dy. 
\end{equation}
We have 
\begin{eqnarray*}
\beta^{\frac{n-1}{n}}\phi'(y)=\left\{\begin{array}{ll} \beta^{\frac{n-1}{n}}\mu^{-1}, & 0<y < \mu M, \\
\alpha_n^{\frac{-1}{n-1}}\beta^{\frac{-1}{n}}  (\cA_{n-1}f(E_{\phi(y)}))^{\frac{1}{n-1}}, & \mu M < y < \|\tilde{\psi}\|_{\infty}. \end{array} \right.
\end{eqnarray*}
Thus, by changing variables with $s= \phi(y)$, and by our choice of $\mu$, 
\begin{eqnarray*}
\int_0^{\infty} (\beta^{\frac{n-1}{n}}\phi'(y))^n\,dy &=&  \int_0^{\mu M}\beta^{n-1} \mu^{-n} \,dy + 
\alpha_n^{\frac{-n}{n-1}}\beta^{-1}\int_{\mu M}^{\|\tilde{\psi}\|_{\infty}} (\cA_{n-1}f(E_{\phi(y)}))^{\frac{n}{n-1}}\,dy \\ 
&=& \beta^{n-1}M \mu^{1-n} + \alpha_n^{-1}\int_M^{\|f\|_{\infty}} \cA_{n-1}f(E_{t}) \,dt \\ 
&\leq& \alpha_n^{-1} \int_0^{\infty} \cA_{n-1}f(E_{t})\,dt \leq 1.  
\end{eqnarray*}
By applying equation (6), page 1080 of \cite{moser:1971iumj} to $\beta^{\frac{n-1}{n}}\phi$, we conclude that \eqref{eq:inparts} is 
bounded from above by a constant depending only on $n$. The proof of Theorem \ref{thm:expint} is complete. 

We finally note that, under the assumptions of Theorem \ref{thm:expint}, the left hand side of \eqref{eq:cava} is finite for every $\beta >0$. 
We fix $M>0$, to be chosen later. After applying Theorem \ref{thm:bemodest} to the right hand term in \eqref{eq:cava}, we need to show that 
\begin{equation}
\label{eq:every}
\int_M^{\infty} s^{\frac{1}{n-1}}\exp \Big(\beta s^{\frac{n}{n-1}}-C\int_M^s\frac{dt}{(\cA_{n-1}f(E_{t}))^{\frac{1}{n-1}}} \Big)\,ds 
\end{equation}
is finite, where $C>0$. 

By H\"older's inequality, 
\begin{equation}
\label{eq:holder}
s-M=\int_M^s \frac{(\cA_{n-1}f(E_{t}))^{\frac{1}{n}}}{(\cA_{n-1}f(E_{t}))^{\frac{1}{n}}}\,dt \leq 
\Big(\int_M^s\frac{dt}{(\cA_{n-1}f(E_{t}))^{\frac{1}{n-1}}}\Big)^{\frac{n-1}{n}} \Big(\int_M^s \cA_{n-1}f(E_{t})\,dt\Big)^{\frac{1}{n}}. 
\end{equation}
By our assumption $\int_0^{\infty} \cA_{n-1}f(E_{t})\,dt$ is finite. Thus, by choosing $M$ large enough so that 
\[
\Big(\int_M^{\infty} \cA_{n-1}f(E_{t})\,dt\Big)^{\frac{-1}{n-1}} > \frac{2\beta}{C}, 
\]
and combining this with \eqref{eq:holder}, we can estimate \eqref{eq:every} from above by 
\[
\int_M^{\infty} s^{\frac{1}{n-1}}\exp(\beta(s^{\frac{n}{n-1}}-2(s-M)^{\frac{n}{n-1}}))\,ds, 
\]
which is clearly finite.



%


%


\def\cprime{$'$}

\end{document}